\documentclass [11pt]{article}
\usepackage{amsfonts}
\usepackage{amsmath}
\newtheorem{theorem}{Theorem}[section]
\newtheorem{proposition}{Proposition}
\newtheorem{lemma}{Lemma}

\begin{document}
\newcommand{\TL}{\tilde{\mathcal{L}^{\mathcal{I}}}}
\newcommand{\LL}{\mathcal{L}^{\mathcal{I}}}
\newcommand{\KK}{\mathcal{K}_{\mathcal{I}}}
\newcommand{\iast}{\ast_{\mathcal{I}}}
\newcommand{\es}{\left( \begin{array}{c} g\ast\phi_s \\ \phi_s\end{array} \right)}
\newcommand{\ei}{\left(\begin{array}{c} g\ast\phi_i \\ \phi_i\end{array}
\right)}
\newcommand{\ees}{\left(\begin{array}{c} g\ast_{\mathcal{I}}\phi_s \\ 0 \end{array}
\right)}
\newcommand{\eei}{\left(\begin{array}{c} g\ast_{\mathcal{I}}\phi_i \\ 0 \end{array}
\right)}
\newcommand{\eees}{\left(\begin{array}{c} 0 \\ \phi_s\end{array}
\right)}
\newcommand{\eeei}{\left(\begin{array}{c} 0 \\ \phi_i\end{array}
\right)}
\topmargin 0pt
\headsep 0pt
\date{}
\title{Poisson Statistics for the Largest Eigenvalues of Wigner Random Matrices with Heavy Tails}
%\footnote{ AMS 2000 subject classification: 
%keywords and phrases: random matrices,  singular% values, Cauchy distribution}

\author{ 
Alexander Soshnikov\thanks{
Department of Mathematics,
University of California at Davis, 
One Shields Ave., Davis, CA 95616, USA.
Email address: soshniko@math.ucdavis.edu.
Research was supported in part by the Sloan Research Fellowship and the 
NSF grant DMS-0103948. }}

\date{}
\maketitle
\begin{abstract}
We study large Wigner random  matrices in the case when
the marginal distributions of matrix entries 
have heavy tails. We prove that the largest eigenvalues of such matrices have Poisson statistics.
\end{abstract}

\section{Introduction and Formulation of Results.}

The main goal of this paper is to study the largest eigenvalues of  Wigner real symmetric and Hermitian 
random matrices in the case when the matrix entries 
have 
heavy tails of distribution.
We remind that a real symmetric Wigner random matrix is defined as a square symmetric
$ n\times n $ matrix with i.i.d. entries up from the diagonal
\begin{equation}
\label{matrica}
 A=(a_{jk}), \ \ a_{jk}=a_{kj}, \ 1 \leq j \leq k \leq n,  \ \ \{a_{jk}\}_{j\leq k} - \ {\rm i.i.d.} \ \ {\rm random} \ \ {\rm variables}.
\end{equation}
A Hermitian Wigner random matrix is defined in a similar way, namely as a square $ n \times n $ Hermitian matrix with i.i.d. entries up
from the diagonal
\begin{equation}
\label{hermitian}
 A=(a_{jk}), \ \ a_{jk}=\overline{a_{kj}}, \ 1 \leq j \leq k \leq n, 
\ \ \{a_{jk}\}_{j\leq k} - \ {\rm i.i.d.} \ \ {\rm complex} \ \ {\rm random} \ \ {\rm variables}
\end{equation}
(it is  
assumed quite often that the distribution of the diagonal matrix entries is different from the distribution of the off-diagonal entries).
While our main results
(with obvious modifications) are valid for the Hermitian random matrices with i.i.d. entries as well, we will restrict ourselves for the 
most of the paper to the real symmetric case.

The ensembles (\ref{matrica}), (\ref{hermitian}) were introduced by E.Wigner in the fifties (\cite{Wig1}, \cite{Wig2}), \cite{Wig3})
who considered the case of centered identically distributed (up from the diagonal) matrix entries with the tail of the distribution decaying 
sufficiently fast so that all moments exist. Wigner proved that under the above conditions 
the mathematical expectation of the empirical distribution function
of $n^{-1/2} \* A $ converges as $ n \to \infty $ to the semicircle law, i.e.
\begin{equation}
\label{semicircle}
\frac{1}{n} \* E \left( \# (\lambda_i \leq x, \ \ 1\leq i \leq n) \right) \to F(x)= \int_{-\infty}^x f(t) \* dt,
\end{equation}
where $\lambda_1,\ldots, \lambda_n \ $ are the eigenvalues of $ n^{-1/2} \* A, \ $
the density of the semicircle law is given by $ f(t)= \frac{1}{\pi \* \sigma^2} \* \sqrt{2 \sigma^2 - x^2}, \ $ 
for $ t \in [-\sqrt{2}\* \sigma, \sqrt{2}\*\sigma] $ and $\sigma^2 $ is the second moment of matrix entries. Wigner's result was 
subsequently strengthened by many mathematicians (see e.g. \cite{MP}, \cite{Ar},  \cite{Pas},  \cite{FK}).  In particular, 
Pastur (\cite{Pas})
and Girko (\cite{Gir}) proved that if $A^{(n)}=(a_{ij}^{(n)}), \ \ 1\leq i,j \leq n, \ $ is an $ n\times n \ $ 
random symmetric matrix with independent (not necessarily identicaly distributed) 
centered entries with
the same second moment $\sigma^2$ , then the necessary and sufficient for the convergence of the empirical distribution function of 
the eigenvalues of $n^{-1/2} \* A^{(n)} $ to the semicircle law has the Lindeberg-Feller form: 
$ \ \ \frac{1}{n^2} \sum_{1\leq i \leq j \leq n} \int_{|x|>\tau \* \sqrt{n}} x^2\* dF_{ij}(x) \to 0, \ \ $  where $ F_{ij}(x) $ is the 
distribution function
of $ a_{ij}^{(n)}.$ 

The simplest case of (\ref{matrica}) from the analytical viewpoint is when the entries 
$ a_{ij}, \ \ 1\leq i \leq j \leq n \ $ are independent  Gaussian $\ N(0, 1+\delta_{ij}) $
random variables. Then the ensemble is called Gaussian Orthogonal Ensemble (GOE). We refer the reader to \cite{M}, chapter 6 for the 
discussion of the GOE ensemble. In a similar fashion the ensemble of Hermitian random matrices with independent and identically distributed 
(up from the diagonal) matrix entries is called  Gaussian Unitary Ensemble (GUE, see \cite{M}, chapter 5) amd the ensemble of Hermitian 
self-dual matrices with quaternion entries is called Gaussian Symplectic Ensemble (GSE, see \cite{M}, chapter 7). In the Gaussian ensembles
one can study in detail the local statistical properties of the eigenvalues both in the bulk and at the edge of the spectrum. In the seminal 
papers  (\cite{TW1}, \cite{TW2})
Tracy and Widom  proved 
that after the proper
rescaling the distribution of the largest eigenvalues in the GOE, GUE and GSE 
ensembles converges to what is now called Tracy-Widom distribution. In particular they proved  that
for the largest eigenvalue of GOE 
\begin{equation}
\label{Tracy}
\lim_{n \to \infty} \Pr \left( \lambda_{max} \leq 2\*\sqrt{n} + \frac{s}{n^{1/6}} \right) \to F_{1}=
\exp\left( -\frac{1}{2} \* \int_s^{+\infty} q(x) + (x-s)\* q^2(x)\* dx \right),
\end{equation}
and for the largest eigenvalue of the GUE ensemble
\begin{equation}
\label{Tracy2}
\lim_{n \to \infty} \Pr \left( \lambda_{max} \leq 2\*\sqrt{n}  + \frac{s}{n^{1/6}} \right) \to F_{2}=
\exp\left( - \int_s^{+\infty}(x-s)\* q^2(x)\* dx \right),
\end{equation}
where $ q(x) \ $ is the solution  of the Painl\'{e}ve II differential equation $ q''(x)=x\*q(x) +2\*q^3(x) \ $ determined by the asymptotics
$ q(x) \sim Ai(x) \ $ at $ x = +\infty.$
It was also established (see \cite{TW1}, \cite{TW2}, \cite{For}) that after the rescaling at the edge of the spectrumas
the $k$-point correlation functions have a limit. In the GUE case the limiting $k$-point correlation functions
has a determinantal form
\begin{equation}
\label{corrgue}
\rho_k(x_1,\ldots, x_k)= \det \left( K(x_i,x_j) \right)_{1 \leq i,j \leq k},
\end{equation}
where
\begin{equation}
\label{airyker}
K(x,y)=K_{Airy}(x,y) = \frac{Ai(x)\*Ai'(y)-Ai'(x)\*Ai(y)}{x-y}
\end{equation}
is a so-called Airy kernel. In the GOE case the limiting $k$-point correlation functions have a pfaffian form and can be written as
\begin{equation}
\label{corrgoe}
\rho_k(x_1,\ldots, x_k)= \left(\det \left( K(x_i,x_j) \right)_{1 \leq i,j \leq k}\right)^{1/2},
\end{equation}
where $K(x,y)$ is a $2\times 2$ matrix kernel such that
\begin{eqnarray}
K_{11}(x,y)&=&K_{22}(y,x)= K_{Airy}(x,y) +\frac{1}{2} \* Ai(x) \*\int_{-\infty}^y Ai(t)\* dt, \\
K_{12}(x,y)&=&-\frac{1}{2} \* Ai(x)\* Ai(y) - \frac{\partial}{\partial y} K_{Airy}(x,y),  \\
K_{21}(x,y)&=& \int_0^{+\infty} \left( \int_{x+u}^{+\infty} Ai(v)\*dv\right) \* Ai(x+u) \* du - \epsilon(x-y) +\frac{1}{2} 
\int_y^x Ai(u) \* du  \nonumber \\
&+& \frac{1}{2} \*\int_x^{+\infty} Ai(u) \* du \* \int_{-\infty}^y Ai(v) \* dv,
\end{eqnarray} 
where $ \epsilon(z)= \frac{1}{2} \* sign(z). \ $
Similar formulas hold for the GSE.

Soshnikov (\cite{So1}, see also \cite{So2}) proved that the Tracy-Widom law for the largest eigenvalues is universal provided the laws of 
the distribution of 
matrix entries are symmetric, all moments exist and do not grow too fast. To the best of our knowledge there are no results proving 
universality in the bulk of the spectrum for real symmetric Wigner matrices. In the Hermitian case Johansson 
(\cite{Jo1}) proved the universality in the bulk under the condition that matrix entries have a Gaussian component.

In this paper we will study the spectral properties of Wigner real symmetric and Hermitian random matrices in the case when the 
matrix entries have heavy tails of distribution. In other words, in addition to the assumption that $ a_{jk} $ are i.i.d. real (complex) 
random variables up from the diagonal ($ 1 \leq j \leq k \leq n $) we also assume that
the probability distribution function satisfies
\begin{equation}
\label{tail}
G(x)=\Pr (|a_{jk}| > x) = \frac{h(x)}{x^{\alpha}}, 
\end{equation}
where
$ 0 < \alpha < 2 $ and $ h(x) $ is a slowly varying function at infinity in a sense of Karamata (\cite{Kar}, \cite{Sen}), in 
other words $h(x)$ is a positive function for all $ x>0,$ such that $ \lim_{x \to \infty} \frac{h(t\*x)}{h(x)}=1 $ for all $ t>0.$
The condition (\ref{tail}) means that the distribution of $\ |a_{ij}| \ $ belongs to the domain of the attraction of a stable distribution
with the exponent $\alpha$ (see e.g.  \cite{IL}, Theorem 2.6.1).
The case $ h(x)= const\*(1+o(1)) $ in (\ref{tail}) was considered on a physical level of rigor by Cizeau and Bouchaud in \cite{CB}, 
who called such set of matrices 
``L\'{e}vy matrices'' (see also \cite{BJJNPZ} and \cite{Janik1} for some physical results on the unitary invariant L\'{e}vy ensembles). 
They argued that the typical eigenvalues of $A$
should be of order of $ n^{\frac{1}{\alpha}}.$  Indeed, such normalization makes the euclidean norm of a typical row of $A$ 
to be (with high probability) of the  
order of constant. Cizeau and Bouchaud suggested a formula for the limiting distribution of the (normalized by 
$ n^{\frac{1}{\alpha}}$) eigenvalues.  The formula for the  limiting spectral distribution has a more complicated density then in
the finite variance case. Namely, they claimed that the density should be given by
\begin{equation}
\label{plotnost}
f(x)= L_{\alpha/2}^{C(x), \beta(x)}(x),
\end{equation}
where $L_{\alpha}^{C,\beta} $ is a density of a centered L\'{e}vy stable distribution  defined through its Fourier 
transform $ \hat{L}(k) :$
\begin{eqnarray}
\label{Levy}
& & L_{\alpha}^{C,\beta}= \frac{1}{2\*\pi}\*\int dk \* \hat{L}(k) \* e^{i\*k\*x}, \\
& & \ln \hat{L}(k)= -C\*|k|^{\alpha}\*\bigl(1 +i\beta\*sgn(k)\*\tan(\pi\*\alpha/2)\bigr),
\end{eqnarray}
and $ C(x), \ \ \beta(x) \ $ satisfy a system of integral equations
\begin{eqnarray}
\label{Cbeta}
& & C(x)=\int_{-\infty}^{+\infty} |y|^{\frac{\alpha}{2}-2}\* L_{\alpha/2}^{C(y), \beta(y)}\bigl(x-\frac{1}{y}\bigr) \* dy, \\
& & \beta(x)= \int_{x}^{+\infty}  L_{\alpha/2}^{C(y), \beta(y)}\bigl(x-\frac{1}{y}\bigr)\* dy.
\end{eqnarray}
Note that the density $f(x) \ $ in (\ref{plotnost}) is not a L\'{e}vy density itself, since $C(x), \ \ \beta(x) \ $ are  
functions of $x.$ It was also argued in \cite{CB} that the density $f(x) \ $ should decay like $  \frac{1}{x^{1+\alpha}} \ $ at
infinity,  which suggests that
the largest eigenvalues of $A$ (in the case (\ref{tail}), $ \ h(x)=const) \ $ should be of order $n^{\frac{2}{\alpha}}, \ $ and not
$ n^{\frac{1}{\alpha}} \ $ as the typical eigenvalues.

 Recently, Soshnikov and Fyodorov (\cite{SF}) used the method of 
determinants to study the largest eigenvalues of a sample covariance matrix $A^t\*A$ in the case when the entries of a rectangular
$m\times n $ matrix $A$ are i.i.d. Cauchy 
random variables. The main result of \cite{SF} states that the largest eigenvalues of $A^t\*A$ are of the order $m^2 \* n^2 $ and
\begin{eqnarray}
\label{sf}
\lim_{n \to \infty} E \left(\det(1 +\frac{z}{m^2 \*n^2} \* A^t \* A )\right)^{-1/2} & =&
\lim_{n \to \infty} E \prod_{i=1}^n (1+ z\* \tilde{\lambda_i})^{-1/2} =
\exp\left(-\frac{2}{\pi} \* \sqrt{z} \right)\\
\label{glavforma}
&= & {\bf E} \prod_{i=1}^{\infty} (1+ z\* x_i)^{-1/2},
\end{eqnarray}
where $\tilde{\lambda_i}, \ i=1, \ldots,n $ are the eigenvalues of $\frac{1}{m^2\*n^2}\* A^t\*A, \ \ z \ $ is a complex number with a 
positive real part, 
the branch of $\sqrt{z}$ on $D=\{ z: \Re z >0 \} $ is such that $ \sqrt{1}=1. \ \ $
$ E $ denotes the mathematical expectation with respect to the random matrix ensemble of sample covariance matrices  and 
${\bf E}$ denotes the mathematical 
expectation with respect to the inhomogeneous Poisson random point process on the positive half-axis with the intensity 
$ \frac{1}{\pi \* x^{3/2}}. $ 
The convergence 
is uniform inside $D$ (i.e. it is unform on  the compact subsets of $D$).

The goal of the paper is to study the distribution of the largest eigenvalues of real symmetric (\ref{matrica}) and Hermitian
(\ref{hermitian})  Wigner random matrices in the case of the heavy tails (\ref{tail}).
Let us introduce the normalization coefficient $b_n $ defined so that
\begin{equation}
\label{bn}
\frac{n^2}{2} \* G(b_n\*x) \to \frac{1}{x^{\alpha}},
\end{equation}
for all positive $x>0, \ $ where $G$ is defined in (\ref{tail}) above.
This can be achieved by selecting $b_n$  to be the infimum of all $t$ for which $ G(t-0) \geq \frac{2}{n^2} \geq G(t+0).$
It is clear that 
\begin{equation}
\label{rostbn}
 n^{\frac{2}{\alpha} -\delta} \ll b_n \ll n^{\frac{2}{\alpha} +\delta} 
\end{equation}
for any $ \delta >0 $ as
$n \to \infty $ and $ \frac{n^2 \* h(b_n)}{b_n^{\alpha}} \to 2 $ as $ n \to \infty.$ 
Let us denote the eigenvalues of $b_n^{-1} \* A$ by $ \lambda_1 \geq \lambda_2 \geq \ldots \lambda_n. $
As we will see below this normalization 
is chosen in such a way that the largest normalized eigenvalues 
are of the order of constant and the vast majority of the eigenvalues go to zero.
The main result of our paper is formulated in the next two theorems.

\begin{theorem}

Let $A$ be a Wigner real symmetric (\ref{matrica}) or Hermitian (\ref{hermitian})
random matrix with the heavy tails (\ref{tail}) and $\lambda_1 $  the largest eigenvalue of $b_n^{-1} \* A.$  Then
\begin{equation}
\label{glav}
\lim_{n \to \infty} \Pr (\lambda_1 \leq x) =  \exp(- \* x^{-\alpha}).
\end{equation}

\end{theorem}

{\bf Remark 1}

It is easy to see (and will be shown later) that the r.h.s. of (\ref{glav}) also gives the limiting distribution of the maximum of the
(properly normalized) 
matrix entries, i.e.

\begin{equation}
\label{maksimum}
\lim_{n \to \infty} \Pr (\max_{1\leq i,j \leq n} b_n^{-1} \* |a_{ij}| \leq  \* x) =  
\exp(- x^{-\alpha}).
\end{equation}

Similarly to Theorem 1.1 one can study the distribution of the second, third, fourth, etc largest eigenvalue of $A$. The following general 
result holds.

\begin{theorem}
Let $A$ be a Wigner real symmetric (\ref{matrica}) or Hermitian (\ref{hermitian})
random matrix  with the heavy tail of the distribution of matrix entries (\ref{tail}). 
Then the random point configuration composed of the positive eigenvalues of $b_n^{-1} \* A$  converges in distribution on the cylinder sets  
to the inhomogeneous Poisson random point process on $ (0, +\infty)  $ with the intensity 
$ \rho(x)= \frac{\alpha}{x^{1+\alpha}}.$
\end{theorem}

{\bf Remark 2} As we will see below the result of the Theorem 1.2 can be reformulated in such a way that for any finite integer $k\geq 1$ 
the joint distribution of the $k $ largest eigenvalues of $b_n^{-1} \* A $ coincides in the limit $ n \to \infty $ with the joint 
distribution of the first $k$ order statistics of the absolute values of the normalized matrix entries 
$ \{ b_n^{-1} \*|a_{ij}|, \ \ 1\leq i \leq j \leq n \}. $ This reformulation is also indicative of the method of proof of the results of 
Theorems 1.1 and 1.2. Of course a similar result about the negative eigenvalues of $b_n^{-1} \* A$ holds true as well.

We remind the reader that a
Poisson random point process on the real line with the locally integrable intensity function $ \rho(x) $ is defined in such a way
that the counting functions (e.g. numbers of particles) for the disjoint intervals $ I_1, \ldots, I_k $ are independent Poison random 
variables with the parameters $ \int_{I_j} \rho(x) \* dx, \ \ j=1, \ldots, k.$ Equivalently, one can define the Poisson random point by 
requiring that the $k$-point correlations functions are given by the products of the one-point correlation functions (intensities), i.e.
$\rho_k(x_1, \ldots, x_k)= \prod_{j=1}^k \rho(x_j).$

{\bf Remark 3} The arguments of the proof are quite soft and can be extended without any difficulty to banded random matrices
and some other classes of Hermitian and real symmetric random matrices with independent entries with heavy tails of distribution.

{\bf Remark 4}  For random Schr\"{o}dinger operators the Poisson 
statistics of the eigenvalues in the localization regime was first proved by Molchanov in \cite{Mo} (see also the paper by Minami 
\cite{Mi}). There is a vast literature on the Poisson statistics of the energy levels of quantum systems in the case of the regular 
underlying dynamics (see e.g. \cite{BT}, \cite{Si}, \cite{Maj}, \cite{Sar}, \cite{CLM}, \cite{Mar1}, \cite{Mar2}, \cite{Mar3}).

The rest of the paper is organized as follows. We prove Theorem 1.1 in section 2. The proof of Theorem 1.2 is 
presented in section 3.

\section{Proof of Theorem 1.1}

We will prove Theorems 1.1 and 1.2 in the real symmetric case. In the Hermitian case the arguments are essentially the same.
We start by considering the distribution of the largest matrix entry. It follows from (\ref{bn}) that
\begin{equation}
\label{maxim}
\Pr (\max_{1\leq i \leq j \leq n} b_n^{-1}\*|a_{ij}| \leq  x) = \left( 1- G(b_n\*x) \right)^{n(n+1)/2} =
\left( 1 - \frac{2}{n^2\* x^{\alpha}} \right)^{n(n+1)/2}\*(1 +o(1)) = \exp(-x^{-\alpha})\*(1+o(1)),
\end{equation}
(see e.g. \cite{LLR}, Theorem 1.5.1).

Let us order the $ N= n(n+1)/2 \ \ $ i.i.d. random variables $ |a_{ij}|, \ \ 1\leq i \leq j \leq n, $ i.e.  
$ |a_{i_1 j_1}| \geq |a_{i_2 j_2}| \geq |a_{i_3 j_3}| \geq \ldots \geq |a_{i_N j_N}|,$  where $(i_l,j_l) $ are the indices of the $l$-th 
largest (in absolute value) matrix entry. We will use the notation $ a^{(l)} = b_n^{-1}\* |a_{i_l j_l}|, \ \ l=1, \ldots, N \ \ $ for the 
normalized $l$-th largest absolute value of matrix entries.

In this paper only the statistical properties 
of a finite number of the largest order statistics will be of importance. In particular, for any finite $k$ the inequalities between
the first $k$ order statistics are strict with probability going to 1 as $ n \to \infty.$
We start with a standard proposition that generalizes (\ref{maxim}).

\begin{proposition}
Let $ 0 < c_1 <d_1 <c_2 <d_2 <\ldots c_k <d_k \leq +\infty, $ and $ I_l=(c_l, d_l), \ \ l=1, \ldots k, $ be disjoint intervals on the 
positive half-line. Then the counting functions 
$ \#(I_l)= \# ( 1 \leq i \leq j \leq n : \ \ b_n^{-1}\*|a_{ij}| \in I_l), \ \ l=1, \ldots, k, $ 
are independent in the limit
$ n \to \infty $ and have Poisson distribution with the parameters $ \mu_l= \int_{I_l} \rho(x) dx, $ where 
$ \rho(x)=\frac{\alpha}{x^{1+\alpha}}.$
\end{proposition}

In other words, for any $\epsilon >0 $ the restriction of the point configuration $ \{ b_n^{-1} \* |a_{ij}| \} $ to the interval
$ [\epsilon, +\infty) $ converges in distribution on the cylinder sets to the inhomogeneous Poisson random point process with the intensity
$\rho(x)=\frac{\alpha}{x^{1+\alpha}}$ (we refer the reader to \cite{DVJ}, sections 2 and 3 for the definition and elementary 
properties of Poisson random point processes). In particular,
\begin{equation}
\label{maximk}
\Pr (a^{(k)} \leq x) \to \exp(-x^{-\alpha}) \* \sum_{l=1}^{k-1} \frac{x^{-l\* \alpha}}{l!}.
\end{equation}

The proof of the proposition is a straightforward generalization of the calculations in
(\ref{maxim})  (see e.g. \cite{LLR} Theorem 2.3.1).

Now we proceed to the proof of Theorem 1.1.  The result of the theorem follows from (\ref{maxim}) if we can show that
$ \frac{\lambda_1}{a^{(1)}} \to 1 $ in probability as $ n \to \infty, $ where as before $ \lambda_1 $ is the largest eigenvalue of
$b_n^{-1} \* A $ and $a^{(1)} = b_n^{-1} \* \max_{1 \leq i \leq j \leq n}  |a_{ij}|.$

We start with an elementary lemma.

\begin{lemma}

a) With probability going to 1 there are no diagonal entries greater (in absolute value) than $b_n^{11/20}$.

b) With probability going to 1 there is no pair $(i,j), \ \ 1 \leq i < j \leq n, $  such that $ |a_{ij}| > b_n^{99/100} $ and
$ |a_{ii}| +|a_{jj}| > b_n^{1/10}.$

c) For any positive constant $\delta >0 $ with probability going to 1 there is no row $ 1 \leq i \leq n $  that has at least two entries 
greater in absolute value than $b_n^{\frac{3}{4} +\delta}.$

\end{lemma}

For the convinience of the reader we sketch the proof of the lemma. We start with part a).
We first remind the reader that by a classical result by Karamata (\cite{Kar}, see also \cite{Sen} for a nice exposition of the subject)
a slowly varying function at infinity can be respresented on the interval $[B, +\infty),$ where $B$ is sufficiently large, as
\begin{equation}
\label{Karamata}
h(x)=\exp\left(\eta(x)+\int_B^{+\infty} \frac{\varepsilon(t)}{t} \* dt \right),
\end{equation}
where $\eta(x) $ is a bounded measurable function which has a limit at infinity and $\varepsilon(x)$ is a continuous function
on $[B, +\infty)$ such that $\lim_{x \to \infty} \varepsilon(x) =0. \ $ Using (\ref{tail}), (\ref{rostbn}) and (\ref{Karamata})
we obtain
\begin{eqnarray}
\label{nerav}
& & \Pr \left(\max_{1\leq i \leq n} |a_{ii}| \leq b_n^{11/20} \right)= \bigl(1 - G(b_n^{11/20})\bigr)^n = 
\left(1 - \frac{h(b_n^{11/20})}{b_n^{\frac{11\*\alpha}{20}}}\right)^n \geq \nonumber \\
&\geq & \exp(-n^{-\frac{1}{40}}) \* (1+o(1)).
\end{eqnarray}

In a similar fashion, in order to prove the statement b) we have to show that
$ \frac{n\*(n-1)}{2} \* G(b_n^{99/100}) \*  G(b_n^{1/10}) \leq n^{-1/20} \*(1+O(1)),$ which again follows from
(\ref{tail}), (\ref{rostbn}) and (\ref{Karamata}).  The proof of part c) is very similar and, therefore, left to the reader.

Let $(i_1,j_1)$ be the indices of the maximal (in absolute value) matrix element. It follows from (\ref{maxim}) and part a) of Lemma 1
that with probability going to 1 one has that $ |a_{i_1 j_1}| \geq b_n^{99/100} $ and $i_1\not=j_1.$ Let $f_1 $ be a unit vector in
$\Bbb R^n $ such that all its coordinates except the $i_1$-th and $j_1$-th are zero, the $i_1$-th coordinate is $1$ and the $j_1$-th 
coordinate is $+1$ if $a_{i_1j_1}$ is nonnegative and $-1$ otherwise. Then one can easily calculate tha value of the quadratic form 
$ \ (A \* f_1, f_1)= |a_{i_1j_1}| + \frac{1}{2}\* a_{i_1i_1} +
\frac{1}{2}\*a_{j_1j_1},$
and it follows from (\ref{maxim}) and Lemma 1 b) that with probability going to 1 the sum of the last two terms is much smaller than the 
first term, and, in particular,
$ (A \* f_1, f_1)= |a_{i_1j_1}|\*(1+o(b_n^{-4/5})).$ Since $f_1$ is a unit vector, we clearly have
$\lambda_1 \geq (A \* f_1, f_1).$
To show that
$|a_{i_1j_1}|\*(1+o(1))$ is also an upper bound on the largest eigenvalue we need one more lemma.

\begin{lemma}

With probability going to 1 a Wigner matrix  (\ref{matrica}), (\ref{tail}) 
does not have a row such that both the maximum of the absolute values of the matrix elements in that 
row and the sum of the absolute values of the remaining elements in the row are at least $b_n^{\frac{1}{2}+\frac{\alpha}{4}}$. In other words,
\begin{equation}
\label{netu}
\Pr \{ \exists i, \ \ 1 \leq i \leq n : \ \ \max_{1\leq j \leq n} |a_{ij}| > b_n^{\frac{3}{4}+\frac{\alpha}{8}}, 
\ \ \left( \sum_{1 \leq j \leq n}|a_{ij}| \right) - 
\max_{1\leq j \leq n} |a_{ij}| >b_n^{\frac{3}{4}+\frac{\alpha}{8}} \} \to 0
\end{equation}
as $ n \to \infty.$
\end{lemma}

We will treat separately the two cases $ 1 < \alpha < 2$ and $ 0 < \alpha \leq 1.$

Let us start with the first case, i.e. $ 1 < \alpha < 2.$  We choose $ \epsilon =\frac{1}{2\*M +1}, $ where $ M $ is some integer, so 
that $\epsilon < \frac{1}{4}- \frac{\alpha}{8}.$ We claim that there exists sufficiently small positive constant $\gamma  \ \ ( \gamma=
\frac{1}{4\*M +2} $ would suffice) that
for any row $ 1 \leq i \leq n $ and for all integers $ 0 \leq k \leq M $ 
\begin{eqnarray}
\label{matozh}
& & E \left(\# \left(1\leq j \leq n :  |a_{ij}| \geq b_n^{k\*\epsilon} \right) \right)= n\* G(b_n^{k\*\epsilon}),  \\
\label{ver}
& & \Pr \left( \# \left(1\leq j \leq n :  |a_{ij}| \geq b_n^{k\*\epsilon} \right) \geq 2 \*  n\* G(b_n^{k\*\epsilon})\right) \leq
\exp(-n^{\gamma})
\end{eqnarray}
for sufficiently large $n.$ Indeed, (\ref{matozh}) immediately follows from the definition of $G$ in (\ref{tail}), while (\ref{ver})
follows from the Chernoff's inequality
$ \Pr( X \geq EX +t) \leq \exp(-\frac{t^2}{2\*n \*p} +\frac{t^3}{6\*(np)^2}) $ for a binomial random variable $X \sim Bin(n,p) $ and
(\ref{tail}), (\ref{rostbn}).
It then follows from (\ref{ver}) that with probability $ 1 - O(n\* \exp(-n^{\gamma})) $ we have that for all rows 
$ 1 \leq i \leq n $ of the matrix $A$ 
\begin{eqnarray}
\label{barsa}
& & \sum_{j: |a_{ij}| \leq b_n^{(M+1)/(2M+1)}} |a_{ij}|  \leq \sum_{k=0}^M  
\# \left(1\leq j \leq n :  |a_{ij}| \geq b_n^{k\*\epsilon} \right) \* b_n^{(k+1)\*\epsilon} \leq \nonumber \\
& & \sum_{k=0}^M b_n^{(k+1)\*\epsilon} \* 2 \*n \* G(b_n^{k\*\epsilon}) =
\sum_{k=0}^M b_n^{(k+1)\*\epsilon} \* 2 \*n \* h(b_n^{k\*\epsilon}) \* b_n^{-\alpha \* k \*\epsilon} \leq 
n \* b_n^{\frac{1}{4}- \frac{\alpha}{8}} \leq  b_n^{\frac{1}{2}+\frac{\alpha}{4}} \leq \frac{1}{2}\*b_n^{\frac{3}{4}+\frac{\alpha}{8}}.
\end{eqnarray}
for sufficiently large $n.$ 

Finally, we observe that for any fixed row $i$
\begin{equation}
\label{odnako1} 
\Pr ( \sum_{j: b_n^{\frac{M+1}{2M+1}} \leq |a_{ij}| \leq b_n^{\frac{3}{4}+\frac{\alpha}{16}}} |a_{ij}| \geq 
\frac{1}{2}\*b_n^{\frac{3}{4}+\frac{\alpha}{8}} ) < 
\exp(-n^{\kappa}),
\end{equation}
for sufficiently small positive $ \kappa >0.$
The proof of (\ref{odnako1}) is similar to the arguments given above and will be left to the reader.

Since by part c) of Lemma 1 with probability going to 1 each row has at most one entry greater in absolute value
than $b_n^{\frac{3}{4}+ \frac{\alpha}{16}} $ the bounds (\ref{barsa}) and (\ref{odnako1} imply the statement of the Lemma 1.2. 
in the case $ 1 < \alpha < 2.$

In the case $ 0 < \alpha < 1 $ the onsideration is similar. 
We again choose $ \epsilon =\frac{1}{2\*M +1}, $ where $ M $ is some integer in such a 
way that $\epsilon < \frac{\alpha}{8} $  and $\gamma= \frac{1}{4\*M +2} $ so that
with probability $ 1 - O(n\* \exp(-n^{\gamma})) $ we have that for all rows 
$ 1 \leq i \leq n $ of the matrix $A$ 
\begin{eqnarray}
\label{barsa1}
& & \sum_{j: |a_{ij}| \leq b_n^{\frac{M+1}{2M+1}}} |a_{ij}|  \leq \sum_{k=0}^M 
\# \left(1\leq j \leq n :  |a_{ij}| \geq b_n^{k\*\epsilon} \right) \* b_n^{(k+1)\*\epsilon} \leq \nonumber \\
& & \sum_{k=0}^M b_n^{(k+1)\*\epsilon} \* 2 \*n \* G(b_n^{k\*\epsilon}) =
\sum_{k=0}^M b_n^{(k+1)\*\epsilon} \* 2 \*n \* h(b_n^{k\*\epsilon}) \* b_n^{-\alpha \* k \*\epsilon} \leq 
b_n^{\epsilon} \* 2 \*n \* b_n^{\frac{1-\alpha}{2}}
\leq  b_n^{\frac{1}{2}+\frac{\alpha}{4}} \leq \frac{1}{2}\*b_n^{\frac{3}{4}+\frac{\alpha}{8}}.
\end{eqnarray}
for sufficiently large $n.$  Combined with (\ref{odnako1}) this finishes the proof.
Lemma 2 is proven.

{\bf Remark 5}
It follows from (\ref{maxim}) and Lemma 2 that 
\begin{equation}
\label{normy}
\|A\|_{\infty}: =\max_{1\leq i \leq n} \sum_{1\leq j \leq n} |a_{ij}|= \max_{1\leq i \leq j \leq n} |a_{ij}| \* (1+o(1))
\end{equation}
with probability going to 1. The Theorem 1.1 now follows from the fact that the matrix norm $ \|A\|_{\infty} $ is an upper bound for the 
largest
eigenvalue of the matrix $A.$

\section{Proof of Theorem 1.2}

Let as in section 2 $ \ \ a^{(l)} = b_n^{-1}\* |a_{i_l j_l}|, \ \ l=1, \ldots, N=n(n+1)/2 \ \ $ be the values $ \ b_n^{-1} \* |a_{ij}|, \ \ 
1\leq i \leq j \leq n \ $ put in the decreasing order. In particular, $ a^{(1)}= \max_{ij} b_n^{-1} \* |a_{ij}|. \ $
Similarly to section 2 we construct  the unit vectors  $f_l  \in
\Bbb R^n, \ \ l=1,2,\ldots, $ such that all coordinates of $f_l$ except the $i_l$-th and $j_l$-th are zero, 
the $i_l$-th coordinate is $1$ and the $j_l$-th 
coordinate is $+1$ if $a_{i_lj_l}$ is nonnegative and $-1$ otherwise. Then $ A \* f_l = |a_{i_lj_l}| \* f_l + 
\frac{1}{\sqrt{2}}\* \sum_{m\not=i,j} (a_{i_l m} +
a_{j_l m}) \* e_m, \ \ $ where $e_m$ are the standard basic vectors in $\Bbb R^n, $ i.e. all coordinates of $e_m $ are zero, 
except for the $m$-th one which is 1. Since $ \left(\sum_{m\not=i,j} a^2_{i_l m} \right)^{1/2} \leq \sum_{m\not=i,j} |a_{i_l m}|, $
it follows from Proposition 1 and Lemma 2 that for any finite number of the largest values $a^{(l)}, \ \ 1\leq l \leq k \ $ one has that
$ b_n^{-1} \* A f_l = a^{(l)} \* f_l  + r_l, $ where the euclidian norms of $r_l$ are o(1).  Therefore by an elementary argument in 
perturbation theory for Hermitian operators 
we obtain that $b_n^{-1} \* A $ has eigenvalues $ a^{(l)}\*(1+o(1)), \ \ 1\leq l \leq k, \ \ $ for any finite $k. \ \ $
The last thing we need to show is that (with probability going to 1) these eigenvalues are exactly the $k$ largest eigenvalues of 
$b_n^{-1}\* A. \ $ In other words, so far we proved that with probability going to 1 one has $\lambda_l \geq a^{(l)}\*(1+o(1)), \ \ 
l=1, \ldots, k.$  Our goal now is to prove the reverse inequalities $\lambda_l \leq a^{(l)}\*(1+o(1)), \ \ 
l=1, \ldots, k.\ $  To simplify the exposition we will prove the result for the second eigenvalue. The reasoning in the general case is very 
similar. Let, as before, $ (i_1,j_1) $ be the indices of the largest in absolute value matrix element. Consider an $(n-1)\times (n-1) $ 
submatrix obtained from $A$ by deleting the $i_1$-th row and the $i_1$-th column. Let us denote the submatrix by $ B.$ 
It follows from Lemma 1 that $i_2 \not= i_1, \ \ 
j_2 \not= j_1 \ $ with probability going to 1. By the same reasoning as in Theorem 1 we have that the largest eigenvalue of $ b_n^{-1} \* B$ 
is $ a^{(2)}(1+o(1)),$ but by interlacing property of the eigenvalues of $A$ and $B$ we have that the largest eigenvalue of $B$ is not 
smaller than the second largest eigenvalue of $A$. Theorem 1.2 is proven.

\noindent

%{\bf Acknowledgements.}\,

\def\cmp{{\it Commun. Math. Phys.} }

\end{document}